\documentclass[11pt]{article}
\usepackage{graphicx,amssymb,mathrsfs,amsmath}
\usepackage{color}
\usepackage{rotating}
\newfont{\bb}{msbm10}

\def\Diag{{\rm Diag}}

\newtheorem{theorem}{Theorem}[section]

\newtheorem{lemma}{Lemma}[section]
\newtheorem{remark}{Remark}[section]
\newtheorem{assumption}{Assumption}[section]

\newcommand{\tr}{^{\sf T}}

\numberwithin{equation}{section}

\begin{document}
\cleardoublepage \pagestyle{myheadings}

\bibliographystyle{plain}

\title{\bf    General parameterized proximal point algorithm with  applications in statistical learning
\thanks{\scriptsize The work was supported by the National Science Foundation of China under
grant 11671318 and the Natural Science Foundation of Fujian province under grant 2016J01028.} \\}

\author
{  Jianchao Bai$^a$\footnote{\scriptsize{Corresponding author.
E-mail address: bjc1987@163.com.} } \quad Jicheng Li$^{a}$  \quad Pingfan Dai$^{a,b}$ \quad Jiaofen Li$^{c}$
 \\ \\
{\small \it $^a$School of Mathematics and Statistics, Xi'an Jiaotong University, Xi'an 710049, P.R.  China}\\
{\small\it  $^{b}$Department of Information Engineering, Sanming University, Sanming  365004, P.R. China}\\
{\small\it   $^c$College of Mathematics and Computational Science,  Guilin University of}\\
{\small\it Electronic Technology, Guilin 541004, P.R. China} \\}

\date{}
\maketitle

\hrule
\bigskip
{\noindent\bf Abstract} \vskip 1mm
\small
In the literature, there are a few researches to design  some parameters in the Proximal Point Algorithm (PPA), especially for the multi-objective convex optimizations. Introducing some parameters to  PPA can make it more  flexible and attractive. Mainly motivated by our recent work (Bai et al., A parameterized proximal point algorithm for separable convex optimization, Optim. Lett. (2017)  doi: 10.1007/s11590-017-1195-9), in this paper we develop a general parameterized
PPA with a relaxation step  for solving the multi-block separable structured convex programming. By making use of the variational inequality and some mathematical identities, the global convergence and  the  worst-case $\mathcal{O}(1/t)$
convergence rate of the proposed algorithm are established. Preliminary numerical experiments on solving a sparse matrix minimization problem from statistical learning validate that our  algorithm is more efficient  than several state-of-the-art algorithms.

\bigskip
\hrule
\bigskip

%==========================================
\section{Introduction}
%==========================================
Throughout this  paper, let $\mathbb{R}(\mathbb{R}^+), \mathbb{R}^m, \mathbb{R}^{m\times n}$ be the set of real (positive) numbers, the set of  $m$ dimensional real column vectors and the set of $m\times n$ dimensional real matrices, respectively. The symbol $\|z\|_2$ denotes the Euclidean norm of $z\in \mathbb{R}^m$, which is defined by $\|z\|_2=\sqrt{\langle z,z\rangle}$ with  the standard inner product $\langle\cdot,\cdot\rangle$. For any symmetric positive definite matrix $G\in\mathbb{R}^{m\times m}$, $\|z\|_G=\sqrt{\langle z,Gz\rangle}$  represents the weighted $G$-norm of $z$. We also use $\tr$ and $\textbf{I}$ to stand  for the transpose of a vector/matrix and  the identity matrix
with proper dimensions, respectively.

Consider the following  multi-block separable structured convex optimization
\begin{equation} \label{Sec1-prob}
\begin{array}{lll}
 \min &\sum\limits_{i=1}^{p}f_i(x_i)\\
\textrm{s.t. } &\sum\limits_{i=1}^{p}A_ix_i =b,\  x_i \in\mathcal{X}_i,
 \end{array}
\end{equation}
where $p> 1$ is a positive integer, $f_i(\cdot):\mathbb{R}^{n_i}\rightarrow\mathbb{R}$ are closed convex functions (possibly nonsmooth);  $A_i\in\mathbb{R}^{m\times n_i}, b\in \mathbb{R}^{m}$ are respectively given matrices and vector; all structured  sets $\mathcal{X}_i\subset \mathbb{R}^{n_i}(i=1,\cdots,p)$ are closed and convex. Throughout the discussions, we make the following  assumption:
\begin{assumption}\label{assumption}
The solution set of the problem (\ref{Sec1-prob}) is  nonempty and all the matrices $A_i(i=1,\cdots,p)$ have full column rank.
\end{assumption}

Note that the first part of  Assumption \ref{assumption} is basic and the second part is necessary. For instance, when dealing with a minimization problem subject to the linearly equality constraint $Ax=b$ involving large-size coefficient matrix and variable, we can split
\[
A=[A_1,\cdots,A_p], \quad   x=(x_1\tr,\cdots,x_p\tr)\tr
\]
to reduce the dimensions of $A$ and $x$,  then in such case each   $A_i$  has full column rank. In fact, many  practical application
problems do contain at least two different  variables (or can be transformed into an equivalent  problem with at least two  variables), e.g. the total-variational image deblurring problems \cite{HeMaYuan2017,ZhangZhu2016},  the transformed joint sparse recovery problem \cite{LiaoYang2015}, the sparse inverse covariance estimation problem \cite{Friedman2008}, the low-rank and sparse problem \cite{LiuLiBai2017,TaoYuan2014} and so forth.

The  Proximal Point Algorithm (PPA), which was originally proposed to tackle the monotone operator inclusion problems \cite{Moreau1965,Martinet1970}, is regarded as a  powerful  benchmark method for solving the convex problems like (\ref{Sec1-prob}). As verified by Rockafellar \cite{Rockafellar1976}, the well-known Augmented Lagrangian Method (ALM) for (\ref{Sec1-prob}) with $p=2$ was actually an application of  PPA to its dual problem. Moreover, the classical Alternating Direction Method of Multipliers (ADMM) can be also treated as a special variant of PPA to the dual problem \cite{EcksteinBertsekas1992}. In the last several years, there has been a constantly increasing interest in  developing the theories and applications of PPA. For example, He et al.\cite{HeYuanZhang2013} showed a customized application of the classical PPA to the model (\ref{Sec1-prob}) with $p=1$, where  some image processing problems were carried out to
show the efficiency of the method therein. Ma and Ni \cite{MaNi2016} proposed a parameterized PPA for (\ref{Sec1-prob}) with $p=1$, where both the basis pursuit problem and the matrix completion problem  were tested in  experiments to exam the numerical performance of their algorithm. Recently, Cai et al.\cite{CaiGuHeYuan2013} designed a PPA with a relaxation step for the model (\ref{Sec1-prob}) with $p=2$, whose global convergence and the worst-case sub-linear convergence rate were analyzed in detail. More recently, by introducing some parameters to the metric proximal matrix, an extended parameterized PPA based on \cite{MaNi2016} was developed for  the  two block separable convex programming \cite{BaiZhang2017}, whose  effectiveness and robustness was demonstrated by testing a sparse vector optimization problem in  the statistical learning  compared with two popular algorithms.

To the best of our knowledge, there are a few researches on the parameterized PPA for solving the multi-block model (\ref{Sec1-prob}) with at least three variables. Based on such observation and  motivated by our recent published  work \cite{BaiZhang2017}, the aims of this article are to design a general  parameterized PPA with a relaxation step (\textbf{GR-PPA})  for solving (\ref{Sec1-prob}) and to test some practical examples having more than two variables to investigate the performance of our GR-PPA. In the remaining parts, Section 2 shows the details  of constructing the proposed GR-PPA
and analyzing its convergence theories. In Section 3,  some numerical experiments are carried out in terms of different parameters,   tolerances and   initial values. Finally, we conclude   the paper in Section 4.

%==========================================
\section{ Main results}
%==========================================
In this section, we first construct a parameterized proximal matrix to design a novel
GR-PPA for (\ref{Sec1-prob}),  where its  convergence is  analyzed in detail afterwards.
The whole convergence analysis is based on the variational inequality and uses some
special techniques to simplify its proof.

%==========================================
\subsection{Formation of GR-PPA }
%==========================================
For any $ \tau\in \mathbb{R}^+$, the Lagrangian function of (\ref{Sec1-prob}) is constructed as
\begin{equation} \label{Sec2-Lagra}
L(x_1,\cdots,x_p,\lambda)= \sum\limits_{i=1}^{p}f_i(x_i)-\tau\left\langle\lambda,
\sum\limits_{i=1}^{p}A_ix_i-b\right\rangle,
\end{equation}
where $\lambda\in \mathbb{R}^m$ denotes the Lagrange multiplier with respect to
the equality constraint. Let $(x_1^*,\cdots,x_p^*,\lambda^*)$ be the saddle-point
belonging to the solution set $\Omega^*$ of (\ref{Sec1-prob}). Then,  the following
basic inequalities
\[
L(x_1^*,\cdots,x_p^*,\lambda)\leq L(x_1^*,\cdots,x_p^*,\lambda^*)\leq L(x_1,\cdots,x_p,\lambda^*),
\]
imply
\[
 \left \{\begin{array}{lll}
 x_1^*=\arg\min\limits_{x_1 \in\mathcal{X}_1} \left\{f_1(x_1)-\tau\left\langle\lambda^*,
 A_1x_1\right\rangle\right\},\\
 \ \vdots\\
 x_p^*=\arg\min\limits_{x_p \in\mathcal{X}_p} \left\{f_p(x_p)-\tau\left\langle\lambda^*,
 A_px_p\right\rangle\right\},\\
\lambda^*=\arg\max\limits_{\lambda\in\mathbb{R}^m} -\tau\left\langle\lambda, \sum\limits_{i=1}^{p}A_ix_i^*-b\right\rangle,
\end{array}\right.
\]
whose optimality conditions are derived as follows
\[
 \left \{\begin{array}{lll}
 x_1^*\in\mathcal{X}_1, & f_1(x_1)-f_1(x_1^*)+\left\langle x_1-x_1^*,
 -\tau A_1\tr\lambda^*\right\rangle\geq 0,
& \forall\ x_1\in\mathcal{X}_1, \\
\quad\vdots &  &  \\
 x_p^*\in\mathcal{X}_p, & f_p(x_p)-f_p(x_p^*)+\left\langle x_p-x_p^*,
 -\tau A_p\tr\lambda^*\right\rangle\geq 0,
& \forall\ x_p\in\mathcal{X}_p,\\
\lambda^*\in\mathbb{R}^m, & \left\langle \lambda-\lambda^*,
\tau(\sum\limits_{i=1}^{p}A_ix_i^*-b)\right\rangle\geq 0,
& \forall\ \lambda\in\mathbb{R}^m.
\end{array}\right.
\]
It is not hard to reformulate the above optimality conditions into a
variational inequality:
\begin{equation}\label{VI}
w^*\in \Omega^*,\quad \textrm{VI}(\phi, \mathcal{J},\Omega):\
\phi(u)-\phi(u^*)+ \left\langle w-w^*,
\mathcal{J}(w^*)\right\rangle \geq 0,\
\forall\ w\in \Omega,
\end{equation}
where
\[
\phi(u)=\sum_{i=1}^{p}f_i(x_i),\quad  \Omega =\mathcal{X}_1\times\cdots\times
\mathcal{X}_p\times \mathbb{R}^m,
\]
\[
u=\left(\begin{array}{c}
x_1\\ x_2\\ \vdots\\  x_p\\
\end{array}\right),\
w=\left(\begin{array}{c}
x_1\\ \vdots\\  x_p\\  \lambda
\end{array}\right) \ \mbox{and} \
\mathcal{J}(w)=\tau\left(\begin{array}{c}
-A_1\tr\lambda\\
\vdots\\
-A_p\tr\lambda\\
\sum\limits_{i=1}^{p}A_ix_i-b
\end{array}\right).
\]
In the following, we also denote
\begin{equation}\label{Sec2-uwJ}
u^k=\left(\begin{array}{c}
x_1^k\\ x_2^k\\ \vdots\\  x_p^k\\
\end{array}\right),\
w^k=\left(\begin{array}{c}
x_1^k\\ \vdots\\  x_p^k\\  \lambda^k
\end{array}\right) \ \mbox{and} \
\mathcal{J}(w^k)=\tau\left(\begin{array}{c}
-A_1\tr\lambda^k\\
\vdots\\
-A_p\tr\lambda^k\\
\sum\limits_{i=1}^{p}A_ix_i^k-b
\end{array}\right).
\end{equation}
Clearly,   we can  get from the skew-symmetric property of
the   mapping $\mathcal{J}(w)$ that
\begin{equation}\label{Sec2-J-Prope}
\left\langle w^{k+1}-w, \mathcal{J}(w^{k+1})\right\rangle= \left\langle w^{k+1}-w,
\mathcal{J}(w)\right\rangle,\ \forall\ w, w^{k+1}\in \Omega.
\end{equation}

It is well-known that the standard PPA  with given iterate $w^k$  reads
the following unified framework:
\begin{equation} \label{Sec2-PPA}
w^{k+1}\in \Omega,\quad \phi(u)- \phi(u^{k+1}) + \left\langle w-w^{k+1},
\mathcal{J}(w^{k+1})+G\left(w^{k+1}-w^k\right)\right\rangle \geq 0,\quad
\forall\ w\in \Omega,
\end{equation}
where  $G$ is  a symmetric positive definite matrix (called the proximal matrix).
Followed by such framework,  a general parameterized proximal matrix  will
be designed for the new GR-PPA.
Concisely, let the matrix
$G$    be of  following block form
\begin{equation} \label{Sec2-G}
G=\left[\begin{array}{cccc|c}
\left(\sigma_1+\frac{\varepsilon^2-1}{s}\right)A_1\tr A_1& & & & -\varepsilon A_1\tr\\
   & \left(\sigma_2+\frac{\tau^2-1}{s}\right)A_2\tr A_2& &  &-\tau A_2\tr\\
   &  &\ddots&  &\vdots\\
   &   && \left(\sigma_p+\frac{\tau^2-1}{s}\right)A_p\tr A_p&-\tau A_p\tr\\ \hline
-\varepsilon A_1 &-\tau A_2& \cdots&-\tau A_p & s\textbf{I}
\end{array}\right],
\end{equation}
where  $(\sigma_1,\cdots,\sigma_p,s)$ are parameters restricted into the domain
\begin{equation} \label{Sec2-para-region}
 \mathcal{K}= \left\{  \sigma_1>\frac{1+(p-1)\tau|\varepsilon|}{s},
 \sigma_i>\frac{1+(p-2)\tau^2+\tau|\varepsilon|}{s},  s>0\  |\  \varepsilon\in
 \mathbb{R},  \tau\in \mathbb{R}^+,  i=2,\cdots,p\right\}.
\end{equation}
For the sake of analysis convenience, here and next we denote
\begin{equation}\label{Sec2-sigma-rho-bar}
\ \bar{\sigma}_i: = \sigma_i+\frac{\tau^2-1}{s},\ \forall\ i=1,\cdots,p.
\end{equation}

\begin{lemma}\label{G-psd}
 Let  $\mathcal{K}$ be defined in (\ref{Sec2-para-region}) and
the matrices $A_i(i=1,2,\cdots,p)$ have full column rank. Then, the matrix $G$ in
(\ref{Sec2-G}) is symmetric positive definite for any  $(\sigma_1,\cdots,\sigma_p,s)\in\mathcal{K}$.
\end{lemma}
\indent{\textbf{Proof} } Notice that the matrix $G$ is symmetric and can be decomposed as
\[
G=D\tr G_0D,
\]
where $D=\Diag(A_1,\cdots,A_p,\textbf{I})$ and
\begin{equation} \label{Sec3-DG0}
  G_0=\left[\begin{array}{cccc|c}
\left(\sigma_1+\frac{\varepsilon^2-1}{s}\right)\textbf{I}& & & & -\varepsilon\textbf{I}\\
   & \left(\sigma_2+\frac{\tau^2-1}{s}\right)\textbf{I}& &  &-\tau\textbf{I}  \\
   &  &\ddots&  &\vdots\\
   &   && \left(\sigma_p+\frac{\tau^2-1}{s}\right)\textbf{I}&-\tau\textbf{I}\\ \hline
-\varepsilon\textbf{I} &-\tau\textbf{I}& \cdots&-\tau\textbf{I} & s\textbf{I}
\end{array}\right].
\end{equation}
Since all the matrices $A_i$ are assumed to have full column rank, therefore,
the matrix $G$ is positive definite if and only if $G_0$ is positive definite.
Using the identity
\begin{equation} \label{Sec2-G1}
 \begin{array}{l}
\left[\begin{array}{ccccc}
\textbf{I} &  & & & \frac{\varepsilon}{s}\textbf{I}\\
 &  \textbf{I}& & & \frac{\tau}{s}\textbf{I}\\
 && \ddots& &\vdots\\
 &&  &\textbf{I}&\frac{\tau}{s}\textbf{I}\\
 & &&  & \textbf{I}
\end{array}\right]  G_0  \left[\begin{array}{ccccc}
\textbf{I} &  & & & \frac{\varepsilon}{s}\textbf{I}\\
 &  \textbf{I}& & & \frac{\tau}{s}\textbf{I}\\
 && \ddots& &\vdots\\
 &&  &\textbf{I}&\frac{\tau}{s}\textbf{I}\\
 & &&  & \textbf{I}
\end{array}\right]\tr\\ \\
 = \left[\begin{array}{ccccc|c}
\left(\sigma_1-\frac{1}{s}\right)\textbf{I} &-\frac{\varepsilon\tau}{s}\textbf{I}&
-\frac{\varepsilon\tau}{s}\textbf{I}&\cdots& -\frac{\varepsilon\tau}{s}\textbf{I}&  \\
-\frac{\tau\varepsilon}{s}\textbf{I} &\left(\sigma_2-\frac{1}{s}\right)\textbf{I}&
-\frac{\tau^2}{s}\textbf{I}&\cdots& -\frac{\tau^2}{s}\textbf{I}&  \\
\vdots&\vdots&\vdots&\ddots&\vdots&\\
-\frac{\tau\varepsilon}{s}\textbf{I}&-\frac{\tau^2}{s}\textbf{I}&-\frac{\tau^2}{s}
\textbf{I}&\cdots& \left(\sigma_p-\frac{1}{s}\right)\textbf{I}& \\ \hline
&& &&& s\textbf{I}
\end{array}\right]
=: \widetilde{G}_0,
\end{array}\end{equation}
it is easy to demonstrate that $G_0$ is positive definite if and only if
$\widetilde{G}_0$ is positive definite, which is guaranteed by the region
$\mathcal{K}$ defined in (\ref{Sec2-para-region}) because in such case both
the upper-left and lower-right block matrices of $\widetilde{G}_0$ are positive
definite.   $\ \   \square$

\begin{remark}
According to the equation (3.7) in \cite{HeXuYuan2016}, the
proximal matrix of their proposed algorithm  can be decomposed as
\[ G=P\tr \bar{G}_0 P  \]
with $P=\Diag(\sqrt{\beta}A_1,\cdots ,\sqrt{\beta}A_p, \mathbf{I}/\sqrt{\beta})$ and
\begin{equation} \label{Sec2-G2}
\bar{G}_0=\left[\begin{array}{cccc|c}
\nu \textbf{I} & -\textbf{I} &\cdots &-\textbf{I} & \textbf{0}\\
-\textbf{I} & \nu \textbf{I} &\cdots &-\textbf{I} & \textbf{0}  \\
\vdots&\vdots&\ddots&\vdots&\vdots\\
-\textbf{I} & -\textbf{I} &\cdots & \nu \textbf{I}& \textbf{0} \\ \hline
\textbf{0} &\textbf{0} &\cdots &\textbf{0} & \textbf{I}
\end{array}\right],
\end{equation}
in which $ \beta>0$
is a penalty parameter with respect  to the equality constraint of (\ref{Sec1-prob})
and $\nu\geq m-1$ denotes the proximal parameter.  Compared
(\ref{Sec2-G1}) to (\ref{Sec2-G2}), if we set
\[
(s,\tau,\varepsilon):=(1,1,1)\quad \&\ \ \sigma_i:=\nu+1,
\quad \forall\ i=1,2,\cdots,p,
\]
then  it is clear that $\widetilde{G}_0=\bar{G}_0,$ which implies that the proximal
matrix of \cite{HeXuYuan2016} in essence is a special case of  (\ref{Sec2-G}).
 That is, the method in \cite{HeXuYuan2016} is  a special case of our GR-PPA without relaxation step.
A similar way can be    used to analyze the proximal matrix in equation (8.6) in
\cite{GuHeYuan2014}. Such relations show that our parameterized proximal matrix is
more flexible and general than some in the past.
\end{remark}

In what follows, the concrete iterative principles of our GR-PPA are analyzed one by one.
Substituting the matrix $G$ into  (\ref{Sec2-PPA}), we obtain
\[
\lambda^{k+1}\in\mathbb{R}^m,\quad \langle \lambda-\lambda^{k+1}, R_{\lambda} \rangle \geq 0,
\quad \forall\ \lambda\in\mathbb{R}^m,
\]
which is equivalent to
\begin{eqnarray*}
0=R_{\lambda}:  =  \tau\left( \sum_{i=1}^{p}A_ix_i^{k+1}-b\right)-\varepsilon A_1(x_1^{k+1}-x_1^k)
-\tau\sum_{i=2}^{p}A_i(x_i^{k+1}-x_i^k)
 +s\left(\lambda^{k+1}-\lambda^{k}\right).
\end{eqnarray*}
Then, it can be derived from the above equality that
\begin{equation} \label{Sec3-labdaupdate}
 \lambda^{k+1}=\lambda^{k}-\frac{1}{s}\left[(\tau-\varepsilon)A_1x_1^{k+1}+
\varepsilon A_1x_1^k+\tau\sum_{i=2}^{p}A_ix_i^k-b\tau  \right].
\end{equation}
 Combining (\ref{Sec2-PPA})-(\ref{Sec2-G}) and (\ref{Sec3-labdaupdate}), we   have
\begin{equation} \label{Sec3-x-subp}
x_1^{k+1}\in\mathcal{X}_1,\quad f_1(x_1)-f_1(x_1^{k+1})+\left \langle
x_1-x_1^{k+1},
R_{x_1}\right\rangle\geq 0, \quad \forall\ x_1\in\mathcal{X}_1,
\end{equation}
where
\[
\begin{array}{lll}
R_{x_1}
&=&  -\tau A_1\tr \lambda^{k+1}-\varepsilon A_1\tr (\lambda^{k+1}-\lambda^{k})
+\left(\sigma_1+\frac{\varepsilon^2-1}{s}\right)A_1\tr A_1(x_1^{k+1}-x_1^k) \\\\
&=& \varepsilon A_1\tr \lambda^{k}-(\tau+\varepsilon)A_1\tr\left\{\lambda^{k}-\frac{1}{s}\left[(\tau-\varepsilon)A_1x_1^{k+1}+
\varepsilon A_1x_1^k+\tau\sum\limits_{i=2}^{p}A_ix_i^k-b\tau  \right]\right\}\\\\
&& +\left(\sigma_1+\frac{\varepsilon^2-1}{s}\right)A_1\tr A_1(x_1^{k+1}-x_1^k)\\\\
&=&-\tau A_1\tr \lambda^{k}+\frac{\tau^2-\varepsilon^2}{s}A_1\tr A_1x_1^{k+1} +\frac{(\tau+\varepsilon)\tau}{s}A_1\tr\left(
\sum\limits_{i=2}^{p}A_ix_i^k-b \right)+\left(\sigma_1+\frac{\varepsilon^2-1}{s}\right)A_1\tr A_1x_1^{k+1}\\\\
&& - \left(\sigma_1+\frac{\varepsilon^2-1}{s}\right)A_1\tr A_1x_1^k+ \frac{(\tau+\varepsilon)\varepsilon}{s}A_1\tr A_1x_1^k\\\\
&=&-\tau A_1\tr \lambda^{k}+\left(\sigma_1+\frac{\tau^2-1}{s}\right)A_1\tr A_1x_1^{k+1}+\frac{(\tau+\varepsilon)\tau}{s}A_1\tr\left(
\sum\limits_{i=2}^{p}A_ix_i^k-b \right)\\\\
&&- \left[\sigma_1+\frac{\tau^2-1-\tau(\tau+\varepsilon)}{s}\right]A_1\tr A_1x_1^k\\\\
&=&-\tau A_1\tr \bar{\lambda}^{k}+\bar{\sigma}_1A_1\tr A_1(x_1^{k+1}-x_1^k)
\end{array}
\]
with $\bar{\sigma}_1$ defined in (\ref{Sec2-sigma-rho-bar}) is positive  by (\ref{Sec2-para-region}) and
\begin{equation} \label{Sec3-x-labmad}
\bar{\lambda}^k=\lambda^{k}-\frac{\tau+\varepsilon}{s}\left(\sum\limits_{i=1}^{p}A_ix_i^{k}-b\right).
\end{equation}
By (\ref{Sec3-x-subp}) together with $R_{x_1}$, obviously,   $x_1^{k+1}$ is the solution of the following  problem
\begin{eqnarray} \label{Sec3-xsubk1}
x_1^{k+1} &=&\arg\min\limits_{x_1\in \mathcal{X}_1}
\left\{
f_1(x_1)-\left\langle A_1x_1,\tau \bar{\lambda}^k\right\rangle + \frac{\bar{\sigma}_1}{2}\left\|A_1(x_1-x_1^k)\right\|^2
\right\} \nonumber \\
&=& \arg\min\limits_{x_1\in \mathcal{X}_1}\left\{f_1(x_1)+ \frac{\bar{\sigma}_1}{2}
\left\|A_1(x_1-x_1^k)-\frac{\tau}{\bar{\sigma}_1}\bar{\lambda}^k \right\|^2\right\}.
\end{eqnarray}
Here we get from (\ref{Sec3-labdaupdate}) and (\ref{Sec3-x-labmad}) that
 \begin{eqnarray}\label{Sec3-x-labmad-1}
\bar{\lambda}^{k+1}
&= &\lambda^{k+1}-\frac{\tau+\varepsilon}{s}\left(\sum\limits_{i=1}^{p}A_ix_i^{k+1}-b\right)  \nonumber \\
& =& \lambda^{k}-\frac{1}{s}\left[(\tau-\varepsilon)A_1x_1^{k+1}+
\varepsilon A_1x_1^k+\tau\sum_{i=2}^{p}A_ix_i^k-b\tau  \right] -\frac{\tau+\varepsilon}{s}
\left(\sum\limits_{i=1}^{p}A_ix_i^{k+1}-b\right) \nonumber \\
& =&  \bar{\lambda}^k + \frac{\tau+\varepsilon}{s}\left(\sum\limits_{i=1}^{p}A_ix_i^{k}-b\right)
- \frac{\tau+\varepsilon}{s}\left(\sum\limits_{i=1}^{p}A_ix_i^{k+1}-b\right) \nonumber \\
&  &  -\frac{1}{s}\left[(\tau-\varepsilon)A_1x_1^{k+1}+
\varepsilon A_1x_1^k+\tau\sum_{i=2}^{p}A_ix_i^k-b\tau +\tau A_1x_1^k-\tau A_1x_1^k \right] \nonumber \\
& = & \bar{\lambda}^{k}-\frac{\tau+\varepsilon}{s}\sum\limits_{i=1}^{p}A_i(x_i^{k+1}-x_i^k)
 -\frac{1}{s}\left[ (\tau-\varepsilon) A_1 (x_1^{k+1} - x_1^k) + \tau \left(\sum_{i=1}^{p}A_ix_i^{k} -b\right) \right].
\end{eqnarray}

Analogously, for $i=2,3,\cdots,p,$ it follows from (\ref{Sec2-PPA})-(\ref{Sec2-G}) and (\ref{Sec3-labdaupdate}) that
\begin{equation} \label{Sec3-y-subp}
x_i^{k+1}\in \mathcal{X}_i,\quad f_i(x_i)-f_i(x_i^{k+1})+\left \langle
x_i-x_i^{k+1},
R_{x_i}\right\rangle\geq 0, \quad \forall\ x_i\in\mathcal{X}_i,
\end{equation}
where
\[
\begin{array}{lll}
R_{x_i}
&=&  -\tau A_i\tr \lambda^{k+1}-\tau A_i\tr (\lambda^{k+1}-\lambda^{k})
+\left(\sigma_i+\frac{\tau^2-1}{s}\right)A_i\tr A_i(x_i^{k+1}-x_i^k) \\\\
&=&\tau A_i\tr \lambda^{k}-2\tau A_i\tr \left\{\lambda^{k}-\frac{1}{s}\left[(\tau-\varepsilon)A_1x_1^{k+1}+
\varepsilon A_1x_1^k+\tau\sum\limits_{i=2}^{p}A_ix_i^k-b\tau  \right] \right\}\\\\
&&+\left(\sigma_i+\frac{\tau^2-1}{s}\right)A_i\tr A_i(x_i^{k+1}-x_i^k)\\\\
&=&-\tau A_i\tr \bar{\lambda}^{k+\frac{1}{2}}+\bar{\sigma}_iA_i\tr A_i(x_i^{k+1}-x_i^k),
\end{array}
\]
and
\begin{eqnarray} \label{lambda-half}
\bar{\lambda}^{k+\frac{1}{2}} &=& \lambda^k-\frac{2}{s}\left[ (\tau-\varepsilon)A_1x_1^{k+1}+
\varepsilon A_1x_1^k+\tau\sum\limits_{i=2}^{p}A_ix_i^k-b\tau \right] \nonumber \\
& =& \bar{\lambda}^k+ \frac{\tau+\varepsilon}{s}\left(\sum\limits_{i=1}^{p}A_ix_i^{k}-b\right)
 -\frac{2}{s}\left[ (\tau-\varepsilon)A_1(x_1^{k+1}-x_1^k)+
 \tau\left(\sum\limits_{i=1}^{p}A_ix_i^k-b\right) \right] \nonumber \\
 & =&\bar{\lambda}^k+ \frac{\tau+\varepsilon-2\tau}{s}\left(\sum\limits_{i=1}^{p}A_ix_i^{k}-b\right)
 - \frac{2(\tau-\varepsilon)}{s}A_1(x_1^{k+1}-x_1^k) \nonumber\\
&=&\bar{\lambda}^k-\frac{\tau-\varepsilon}{s} \left[2A_1(x_1^{k+1}-x_1^k)+\sum\limits_{i=1}^{p}A_ix_i^k-b \right].
\end{eqnarray}
Therefore,  $x_i^{k+1}$ is the solution of the following   problem
\begin{eqnarray}\label{Sec3-ysubk1}
x_i^{k+1} &= &\arg\min\limits_{x_i\in \mathcal{X}_i}
\left\{ f_i(x_i)-\left\langle A_ix_i,\tau\bar{\lambda}^{k+\frac{1}{2}}\right\rangle
 + \frac{\bar{\sigma}_i}{2} \left\|A_i(x_i-x_i^k)\right\|^2
  \right\} \nonumber \\
& = & \arg\min\limits_{x_i\in \mathcal{X}_i}\left\{f_i(x_i)+ \frac{\bar{\sigma}_i}{2}
\left\|A_i(x_i-x_i^k)-\frac{\tau}{\bar{\sigma}_i}\bar{\lambda}^{k+\frac{1}{2}} \right\|^2\right\},\ i=2,\cdots,p.
\end{eqnarray}

Consequently, the algorithmic framework of   GR-PPA is described in Algorithm 2.1, where we use $\widetilde{w}^k=(\widetilde{x}_1^k,\cdots,\widetilde{x}_p^k,\widetilde{\lambda}^k)$ to replace the original output of the $x_i$-subproblems in (\ref{Sec3-xsubk1}), (\ref{Sec3-ysubk1}) and Lagrange multipliers with given iterate $w^k=(x_1^k,\cdots,x_p^k,\bar{\lambda}^k)$, and  we use $w^{k+1}=(x_1^{k+1},\cdots,x_p^{k+1},\bar{\lambda}^{k+1})$ to stand for the new iterate after adding a relaxation step.

\begin{remark}
Noticing that the  steps 4-9 in Algorithm 2.1 are actually the PPA updates whose proximal term is implicitly simplified into the quadratic term  of each subproblem.   Algorithm 2.1 is an extension of our  method \cite{BaiZhang2017} for  the  two-block separable  convex problem, but the domain of the parameters restricted into (\ref{Sec2-para-region}) is not a direct extension of the past. Besides, Algorithm 2.1 adopts a serial idea between $x_1$-subproblem and other subproblems, while the parallel idea is used among the  $x_i$-subproblems ($i=2,\cdots,p$).  From the relaxation step of Algorithm 2.1, we immediately have the following certain relationship
\begin{eqnarray}\label{Sec2-121}
w^{k+1}-w^k=\gamma(\widetilde{w}^k-w^k).
\end{eqnarray}
\end{remark}

\vskip2mm
\hrule\vskip2mm
\noindent {\bf Algorithm 2.1} (GR-PPA for solving Problem (\ref{Sec1-prob}))
\vskip1.5mm\hrule\vskip2mm
\noindent 1\ \ \ Choose parameters $(\sigma_1,\cdots,\sigma_p,s)\in \mathcal{K}, \gamma\in(0,2)$
and initialize $(x_1^0,\cdots,x_p^0, \lambda^0)\in \Omega$;\\
2\ \ \   $r^0 = \sum\limits_{i=1}^{p}A_ix_i^0-b;\ \bar{\lambda}^0 = \lambda^{0}
-\frac{\tau+\varepsilon}{s} r^0$ by (\ref{Sec3-x-labmad});\\
3\ \ \ \textbf{For} $k=0,1,\cdots,$  if not converge,  \textbf{do}\\
4\quad \qquad $\widetilde{x}_1^k=\arg\min\limits_{x_1\in \mathcal{X}_1}\left\{f_1(x_1)+ \frac{\bar{\sigma}_1}{2}
\left\|A_1(x_1-x_1^k)-\frac{\tau}{\bar{\sigma}_1} \bar{\lambda}^{k} \right\|^2\right\}$;\\
5\quad \qquad  $r^{k} = \sum\limits_{i=1}^{p}A_ix_i^k-b;\ \Delta x_1^k = \widetilde{x}_1^k - x_1^k$;\\
6\quad \qquad $\bar{\lambda}^{k+\frac{1}{2}}=\bar{\lambda}^k-\frac{\tau-\varepsilon}{s}
\left(2A_1\Delta x_1^k+r^{k} \right)$;\\
7\quad \qquad Update the $x_i$-subproblem for $i=2,\cdots,p$ by (\ref{Sec3-ysubk1}):
\[
\left \{\begin{array}{lll}
 \widetilde{x}_2^k=\arg\min\limits_{x_2\in \mathcal{X}_2}\left\{f_2(x_2)+ \frac{\bar{\sigma}_2}{2}
\left\|A_2(x_2-x_2^k)-\frac{\tau}{\bar{\sigma}_2} \bar{\lambda}^{k+\frac{1}{2}} \right\|^2\right\};\\
\   \vdots\\
\widetilde{x}_p^k=\arg\min\limits_{x_p\in \mathcal{X}_p}\left\{f_p(x_p)+ \frac{\bar{\sigma}_p}{2}
\left\|A_p(x_p-x_p^k)-\frac{\tau}{\bar{\sigma}_p} \bar{\lambda}^{k+\frac{1}{2}} \right\|^2\right\};
\end{array}\right.
\]\\
8\quad \qquad  $\Delta x_i^k = \widetilde{x}_i^k - x_i^k,\ \forall\ i=2,\cdots,p$; \\
9\quad \qquad  $\widetilde{\lambda}^k=\bar{\lambda}^{k}-\frac{\tau+\varepsilon}{s} \sum\limits_{i=1}^{p} A_i\Delta x_i^k
-\frac{1}{s}\left[ (\tau-\varepsilon) A_1 \Delta x_1^k +  r^k \tau\right]$  by (\ref{Sec3-x-labmad-1});\\
10\ \ \quad\quad Relaxation step:
\[
\left(\begin{array}{c}
x_1^{k+1}\\\vdots\\  x_p^{k+1}\\
\end{array}\right)=
\left(\begin{array}{c}
x_1^{k}\\\vdots\\  x_p^{k}\\
\end{array}\right) + \gamma\left(\begin{array}{c}
\Delta x_1^k \\\vdots\\ \Delta x_p^k\\
\end{array}\right)\ \textrm{and}\  \bar{\lambda}^{k+1} =  \bar{\lambda}^{k}
+ \gamma (\widetilde{\lambda}^k - \bar{\lambda}^{k}).
\]
\hrule\vskip4mm

%==========================================
\subsection{Convergence analysis of GR-PPA }
%==========================================
This subsection analyzes the global convergence and the ergodic convergence rate of Algorithm 2.1. First of all, we present  an important lemma described as follows.
\begin{lemma}\label{Contract}
The sequence $\{w^{k+1}\}$ generated by Algorithm 2.1 satisfies
\begin{equation} \label{Sec3-IEq}
\left\|w^{k+1}-w^*\right\|_G^2\leq \left\|w^{k}-w^*\right\|_G^2- \frac{2-\gamma}{\gamma}
\left\|w^k-w^{k+1}\right\|_G^2,\quad  \forall\ w^*\in \Omega^*.
\end{equation}
\end{lemma}
\indent{\textbf{Proof} } By Algorithm 2.1 and (\ref{Sec2-PPA}), we have
\begin{equation} \label{Sec3-IEqxx}
\widetilde{w}^k\in \Omega,\quad \phi(u)- \phi(\widetilde{u}^{k}) + \left\langle w-\widetilde{w}^k,
\mathcal{J}(\widetilde{w}^k)+G\left(\widetilde{w}^k-w^k\right)\right\rangle \geq 0,\quad
\forall\ w\in \Omega,
\end{equation}
which, by using (\ref{VI}) and (\ref{Sec2-J-Prope}) with setting $w=w^*$, leads to
\begin{equation} \label{Sec3-basic}
\left\langle \widetilde{w}^k-w^*, G\left(w^k-\widetilde{w}^k\right)\right\rangle \geq 0.
\end{equation}
Then, it follows from (\ref{Sec2-121}) and (\ref{Sec3-basic}) that
\begin{eqnarray}
\left\|w^{k}-w^*\right\|_G^2- \left\|w^{k+1}-w^*\right\|_G^2&=& \left\|w^{k}-w^*\right\|_G^2
- \left\|w^k-w^*+w^{k+1}-w^k\right\|_G^2  \nonumber \\
& =&   2\gamma \left\langle w^k-w^*, G(w^k-\widetilde{w}^k)\right\rangle-\gamma^2\left\|\widetilde{w}^k-w^k\right\|_G^2  \nonumber \\
&=&2\gamma \left\langle w^k-\widetilde{w}^k+\widetilde{w}^k-w^*, G(w^k-\widetilde{w}^k)\right\rangle
-\gamma^2\left\|w^k-\widetilde{w}^k\right\|_G^2  \nonumber\\
&=& \gamma(2-\gamma)\left\|w^k-\widetilde{w}^k\right\|_G^2+ 2\gamma \left\langle
\widetilde{w}^k-w^*, G(w^k-\widetilde{w}^k)\right\rangle\nonumber\\
&\geq& \gamma(2-\gamma)\left\|w^k-\widetilde{w}^k\right\|_G^2
= \frac{2-\gamma}{\gamma}\left\|w^k-w^{k+1}\right\|_G^2,\nonumber
\end{eqnarray}
which immediately completes the proof. $\ \   \square$

Lemma \ref{Contract} shows that the sequence $\{w^{k+1}-w^*\}$ is strictly contractive under the weighted $G$-norm. Moreover, the following global convergence theorem holds.
\begin{theorem}\label{global-converge}
Let the parameters $(\sigma_1,\cdots,\sigma_p,s)\in\mathcal{K}$ and the sequence  $\{w^{k+1}\}$  be generated by  Algorithm 2.1. Then, under   Assumption \ref{assumption}  there exists a point $w^\infty \in \Omega^*$ such that
\begin{equation} \label{point-converge}
\lim_{k \to \infty} w^{k+1} = w^\infty.
\end{equation}
\end{theorem}
\indent{\textbf{Proof} }
See the proof of Theorem 1 \cite{BaiZhang2017}.  $\ \   \square$

In order to establish the convergence rate of Algorithm 2.1 in an ergodic sense, we first need to characterize the solution set of $\textrm{VI}(\phi, \mathcal{J},\Omega)$, which had been given by e.g. \cite{HeMaYuan2017} in the following:
\begin{theorem} \label{opt-thm}
The solution set of  $\textrm{VI}(\phi, \mathcal{J},\Omega)$ in (\ref{VI}) is convex and can be characterized as
\[
\Omega^*=\bigcap_{w\in \Omega}\left\{\widehat{w}\in \Omega|\ \phi(u)-
\phi(\widehat{u})+ \left\langle w-\widehat{w}, \mathcal{J}(w)\right\rangle\geq 0\right\}.
\]
\end{theorem}

\begin{theorem}\label{conv-rate} For any $(\sigma_1,\cdots,\sigma_p,s)\in\mathcal{K}$,
let
\[
\textbf{w}_t: =\frac{1}{1+t}\sum_{k=0}^{t}\widetilde{w}^{k}\quad \mbox{and} \quad \textbf{u}_t: =\frac{1}{1+t}\sum_{k=0}^{t}\widetilde{u}^{k},
\]
where  $\{\widetilde{w}^{k}\}$ is the iterative sequence of Algorithm 2.1. Then, under   Assumption \ref{assumption}  we have
\begin{equation} \label{erg-rate}
\phi(\textbf{u}_t)-\phi(u)+ \left\langle \textbf{w}_t-w, \mathcal{J}(w)\right\rangle
\leq\frac{1}{2\gamma(1+t)}\left\|w^{0}-w\right\|_G^2,
\ \forall\ w\in \Omega.
\end{equation}
\end{theorem}
\indent{\textbf{Proof} }
Combining (\ref{Sec3-IEqxx}) with (\ref{Sec2-J-Prope}), we   get
\begin{equation}\label{Sec2-FianlVI}
\phi(u)- \phi(\widetilde{u}^{k}) + \left\langle w-\widetilde{w}^k, \mathcal{J}(w)\right\rangle \geq \left\langle \widetilde{w}^k-w, G\left(\widetilde{w}^k-w^k\right)\right\rangle.
\end{equation}
Meanwhile, by (\ref{Sec2-121}) and the following identity
\[
2\langle a-b, G(c-d)\rangle=\|a-d\|_G^2-\|a-c\|_G^2+\|c-b\|_G^2- \|d-b\|_G^2
\]
with substitutions $(\widetilde{w}^k,w,w^{k+1},w^k)=(a,b,c,d)$,  it follows  that
\begin{equation}\label{Sec22-1}
 \begin{array}{lll}
\left\langle \widetilde{w}^k-w, G\left(\widetilde{w}^k-w^k\right)\right\rangle&=&
\frac{1}{\gamma}\left\langle \widetilde{w}^k-w, G\left(w^{k+1}-w^k\right)\right\rangle
\\
&=&\frac{1}{2\gamma}\left( \left\|\widetilde{w}^k-w^k\right\|_G^2-\left\|\widetilde{w}^k-w^{k+1}\right\|_G^2+ \left\|w^{k+1}-w\right\|_G^2
- \left\|w^{k}-w\right\|_G^2\right),
\end{array}
\end{equation}
where the first two terms
\begin{equation}\label{Sec22-2}
\begin{array}{lll}
\left\|\widetilde{w}^k-w^k\right\|_G^2-\left\|\widetilde{w}^k-w^{k+1}\right\|_G^2&=&
\left\|\widetilde{w}^k-w^k\right\|_G^2-\left\|w^k-\widetilde{w}^k+w^{k+1}-w^k\right\|_G^2\\\\
&=& \left\|\widetilde{w}^k-w^k\right\|_G^2-\left\|w^k-\widetilde{w}^k+\gamma(\widetilde{w}^k-w^k)\right\|_G^2\\\\
&=&\gamma(2-\gamma)\left\|(\widetilde{w}^k-w^k)\right\|_G^2\geq 0.
\end{array}
\end{equation}
Combining (\ref{Sec2-FianlVI})-(\ref{Sec22-2}), we deduce
\[
\phi(u)-\phi(\widetilde{u}^{k})+ \left\langle w-\widetilde{w}^k, \mathcal{J}(w)\right\rangle +\frac{1}{2\gamma}\left\|w^{k}-w\right\|_G^2
\geq \frac{1}{2\gamma}\left\|w^{k+1}-w\right\|_G^2.
\]
Summing the above inequality over $k=0,1,\cdots,t$, we obtain
\[
(1+t)\phi(u)-\sum\limits_{k=0}^{t}\phi(\widetilde{u}^{k})+
\left\langle (1+t)w-\sum\limits_{k=0}^{t}\widetilde{w}^k, \mathcal{J}(w)\right\rangle
+\frac{1}{2\gamma}\left\|w^{0}-w\right\|_G^2 \geq 0,
\]
which by the definitions of $\textbf{w}_t$ and $\textbf{u}_t$ results in
\begin{equation}\label{Sec2-FianlEq}
\frac{1}{1+t}\sum\limits_{k=0}^{t}\phi(\widetilde{u}^{k})-\phi(u)+ \left\langle \textbf{w}_t-w, \mathcal{J}(w)\right\rangle
\leq\frac{1}{2\gamma(1+t)}\left\|w^{0}-w\right\|_G^2.
\end{equation}
Because of the convexity of the function $\phi(u)$ (since all $f_i$ are assumed to be convex), the following inequality holds
\[
\phi(\textbf{u}_t)\leq \frac{1}{1+t}\sum_{k=0}^{t}\phi(\widetilde{u}^{k}),
\]
which, by substituting it into (\ref{Sec2-FianlEq}), completes the whole proof.   $\ \  \square$

\begin{remark}
Theorem \ref{conv-rate} illustrates the worst-case $\mathcal{O}(1/t)$
convergence rate of Algorithm 2.1  in an ergodic sense. By the  region $\gamma\in(0,2)$ in Algorithm 2.1
and the inequality (\ref{erg-rate}),   one may choose a larger value $\gamma$ approximating to 2 so
that the value of right-hand  of (\ref{erg-rate}) is much smaller.
\end{remark}

\section{Numerical experiments}
In this section, we investigate  the  numerical performance of
our proposed GR-PPA for solving a popular sparse matrix optimization problem.
All experiments are  simulated in MATLAB 7.14 (R2012a) on a lenovo-PC with
Intel Core i5 processor (2.70GHz) and 4 GB memory.

\subsection{Test problem}
The   Latent Variable Gaussian Graphical Model Selection (LVGGMS) problem \cite{Chandrasekaran2012,Ma2017} arising in   statistical learning
 is of the following form:
\begin{equation}\label{Sec3-prob}
\begin{array}{lll}
\min  & F(X,S,L):=\langle X, C\rangle-\log\det(X)+ \nu\|S\|_1+\mu tr(L) \\
 \textrm{s.t. } &   X-S+L=\textbf{0},\ L\succeq\textbf{0},
\end{array}
\end{equation}
where $C\in\mathbb{R}^{n\times n}$ is  a given covariance matrix obtained from  the  sample variables,
$\nu$ and $\mu$ are  given positive weighting  factors, $tr(L)$ denotes the trace of the matrix $L$,
$\|S\|_1=\sum_{i,j=1}^n|S_{ij}|$ stands for the $l_1$-norm of the matrix $S$ and $S_{ij}$ means its  $ij$-th entry.

Clearly, the LVGGMS problem (\ref{Sec3-prob}) can be regarded as a special case of (\ref{Sec1-prob}). And
 applying Algorithm 2.1 it is easy to write the three corresponding subproblems   as the following
\begin{equation}\label{Sec3-reformu}
\left \{\begin{array}{lll}
 \widetilde{X}^{k}=\arg\min\limits_{X\in \mathbb{R}^{n\times n}}\left\{\langle X, C\rangle-\log\det(X)+ \frac{\bar{\sigma}_1}{2}
\left\|X-(X^k+\frac{\tau}{\bar{\sigma}_1} \bar{\lambda}^{k}) \right\|_F^2\right\},\\
\widetilde{S}^k=\arg\min\limits_{S\in \mathbb{R}^{n\times n}}\left\{\nu\|S\|_1+ \frac{\bar{\sigma}_2}{2}
\left\|S-(S^k-\frac{\tau}{\bar{\sigma}_2} \bar{\lambda}^{k+\frac{1}{2}}) \right\|_F^2\right\},\\
\widetilde{L}^k=\arg\min\limits_{L\succeq\textbf{0}}\left\{\mu tr(L)+ \frac{\bar{\sigma}_3}{2}
\left\|L-(L^k+\frac{\tau}{\bar{\sigma}_3} \bar{\lambda}^{k+\frac{1}{2}}) \right\|_F^2\right\}.
\end{array}\right.
\end{equation}
Observe that  the above subproblems  have closed formula solutions.
According to the first-order optimality condition of the
$X$-subproblem in  (\ref{Sec3-reformu}), we   derive
\[\begin{array}{l}
\textbf{0}= C-X^{-1}+\bar{\sigma}_1\left(X-X^k-\frac{\tau}{\bar{\sigma}_1} \bar{\lambda}^{k}\right)\\
\Longleftrightarrow\bar{\sigma}_1 X^2 +\left(C-\bar{\sigma}_1X^k-\tau\bar{\lambda}^{k}\right)X-\textbf{I}=\textbf{0}.
\end{array}\]
Then, by using the eigenvalue decomposition
\[
U\textrm{Diag}(\rho)U\tr = C-\bar{\sigma}_1X^k-\tau\bar{\lambda}^{k},
\]
where $\textrm{Diag}(\rho)$ is a diagonal matrix with diagonal entries
$\rho_i (i=1, \cdots, n)$, we get its explicit solution
\begin{equation}
\widetilde{X}^{k}=U\textrm{Dia}g(\gamma)U \tr,
\end{equation}
in which  $\textrm{Diag}(\gamma)$ is the diagonal matrix with diagonal entries
\[
\gamma_i=\frac{-\rho_i+\sqrt{\rho_i^2+4\bar{\sigma}_1}}{2\bar{\sigma}_1}, \quad i=1,2,\cdots,n.
\]
Applying  the soft shrinkage operator $\textrm{Shrink}(\cdot,\cdot)$, see e.g.\cite{TaoYuan2014},
the solution of  the $S$-subproblem is
\begin{equation}
\widetilde{S}^{k}= \textrm{Shrink}\left(S^k-\frac{\tau}{\bar{\sigma}_2}
\bar{\lambda}^{k+\frac{1}{2}}, \frac{\nu}{\bar{\sigma}_2}\right).
\end{equation}
Besides,  it is obvious that the  $L$-subproblem in (\ref{Sec3-reformu}) is equivalent to
\begin{equation}
L^{k+1}=\arg\min\limits_{L\succeq \textbf{0}}
\frac{\bar{\sigma}_3}{2}\left\|L-\widetilde{L}\right\|_F^2
= V\Diag(\max\{\widetilde{\rho}, \textbf{0}\})V\tr,
\end{equation}
where $\max\{\widetilde{\rho}, \textbf{0}\}$ is taken component-wise and
$V\textrm{Diag}(\widetilde{\rho})V\tr$ is the eigenvalue decomposition of the matrix
\[
\widetilde{L}=L^k+\frac{\tau  \bar{\lambda}^{k+\frac{1}{2}}- \mu \textbf{I}}{\bar{\sigma}_3}.
\]

\subsection{Numerical results}
In the following experiments, the parameters $(\nu,\mu)=(0.005,0.05)$ in (\ref{Sec3-prob}) come from \cite{Ma2017} and
the   matrix $C$ is randomly generated by the MATLAB codes of Boyd's homepage\footnote
{http://web.stanford.edu/$\sim$boyd/papers/admm/covsel/covsel\_example.html.}
with $n=100$. The maximal number of
iterations is set as $1000$   and the  following termination criterions are simultaneously used for all algorithms:
\begin{eqnarray*}
\textrm{IER(k)} &:= & \max\left\{\frac{\left\|X^k-X^{k-1}\right\|_F}{\left\|X^k\right\|_F},
\frac{\left\|S^k-S^{k-1}\right\|_F}{\left\|S^k\right\|_F},  \frac{\left\|L^k-L^{k-1}\right\|_F}{\left\|L^k\right\|_F}\right\}\leq \epsilon_1, \\
\textrm{OER(k)} &:= & \frac{|F(X^{k},S^{k},L^{k})-F^*|}{|F^*|}\leq \epsilon_2,\\
\textrm{CER(k)} &:= &\frac{\left\|X^{k}-S^{k}+L^{k}\right\|_F}{\max\left\{1, \left\|X^k\right\|_F,
\left\|S^k\right\|_F, \left\|L^k\right\|_F\right\}}\leq\epsilon_3,
\end{eqnarray*}
where  $(\epsilon_1,\epsilon_2,\epsilon_3)$ are  given tolerances, $(X^k,S^k,L^k)$ are the $k$-th iterative values and $F^*$ is the approximate optimal objective
function value by running Algorithm 2.1 after 1000 iterations.  The   penalty parameter and the relaxation factor in all involved algorithms are set as $(\beta, \gamma) = (0.05,1.8).$ We   choose the feasible points $(X^0,S^0, L^0,\lambda^0)=(\textbf{I},\textbf{4I},3\textbf{I},\textbf{0})$ as   the initial  starting iterates.

First of all, for the two-block classical {\it lasso} problem, as shown in \cite{BaiZhang2017} that  Algorithm 2.1 could perform significantly better than   ADMM and  the customized relaxed PPA \cite{GuHeYuan2014}  when properly choosing the parameters, especially for solving large-size problems and high accurate solutions are required.   Followed by the similar way of adjusting the  parameters \cite{BaiZhang2017}, we first would like to  investigate the effects of the parameters $(\sigma_1,\sigma_2,\sigma_3,s)$ on Algorithm 2.1 for solving the three-block case (\ref{Sec3-prob}). For this purpose, given $(\varepsilon, \tau)$ we randomly choose some values of one parameter  and fix the remaining  as known values, then it is clear from (\ref{Sec2-para-region}) that
\begin{equation} \label{Sec3-region}
  s>0,\ \sigma_1>\frac{1+2\tau|\varepsilon|}{s},\
 \sigma_i>\frac{1+\tau^2+\tau|\varepsilon|}{s},\ \forall i=2,3.
\end{equation}
Here we can note that  $\sigma_1$ and $\sigma_i(i=2,3)$ play  the same role  if taking $\varepsilon=\tau$, since in such case their regions are the same. Therefore, in the following  we would take $\varepsilon=\tau=\frac{\sqrt{5}-1}{2}$ (the golden section ratio) to adjust the relative better value of e.g. $\sigma_1 $ for an example.  We also investigate the effect of the  parameter $s$ on Algorithm 2.1.

Under the tolerance $\epsilon_i=10^{-8}\ (i=1,2,3)$, experimental results of adjusting the parameters $\sigma_1$ and $s$ are respectively reported in Tables 1-2, where the parameter   in Table 1 is  restricted by $\sigma_1>\frac{2+(\sqrt{5}-1)^2}{20}\approx  0.1764$  and the parameter in Table 2 is restricted by $s>\frac{2+(\sqrt{5}-1)^2}{0.178}\approx   9.9097$. The notations ``IT", ``CPU"  denote the number of iterations and the CPU time in seconds, respectively. We can observe from Tables 1-2 that:
\begin{itemize}
\item
For the parameters $(\sigma_1,\sigma_2,\sigma_3, s)$, the  experimental results in each column of IER  and CER are nearly the same when fixing any three parameters with one parameter changing;
\item
With the increase of the parameter $\sigma_1$ ($s$), both the
iteration number and the CPU time tend to increase, which verifies the comments in the final part of  \cite{HeXuYuan2016}. That is, slow convergence would occur in terms of solving the  proximal subproblem if the proximal parameter is set too large, which  may significantly
affect  the whole computational efficiency of the algorithm;
\item   Both Tables 1-2 show that
$(\sigma_1,\sigma_2,\sigma_3,s)=\left(0.178,0.178,0.178, 10\right)$
would
be  a  relative reasonable   choice for tackling the test problem    (\ref{Sec3-prob}) and we would use them as the default parameter setting for GR-PPA.
\end{itemize}
\vskip2mm
\begin{center}
\begin{tabular}{ccccccccccc}
\hline
 $\sigma_1$&&  IT&&  CPU &&IER   &&OER &&CER  \\
\hline
 0.178 && 160  && 4.82  && 4.30e-9 &&1.27e-12  && 6.97e-9\\
 0.200 && 163  &&  5.24   && 2.73e-9  &&1.69e-12  && 5.71e-9\\
 0.350 && 168  && 5.44  &&  2.43e-9&&  2.63e-12&&  7.75e-9  \\
 0.500 && 170 && 5.65 &&  5.65e-9&& 1.15e-11&&  9.93e-9 \\
  0.600 && 182 &&  5.75 &&  3.06e-9&&  7.13e-11&&  8.58e-9 \\
 0.800 &&  217 &&  6.53  &&  1.85e-9&& 1.39e-10&& 7.19e-9\\
   0.900 && 236 &&  7.36 &&  1.09e-9&&  1.62e-10&&  7.21e-9 \\
 1.000 && 244  &&   7.43 &&  1.81e-9&& 1.93e-10&&  7.49e-9 \\
   3.000 &&318  && 10.01 &&  9.72e-9&&  2.06e-10&&  8.07e-9 \\
   5.000 &&505  && 15.43 &&  9.82e-9&&  2.85e-12&&  4.55e-9 \\
   8.000 &&775  &&  23.98 &&  9.85e-9&&  2.80e-12&&  4.56e-9 \\
   10.000 &&950  &&  30.44 &&  9.88e-9&&  2.40e-12&&  4.56e-9 \\
\hline
\end{tabular}
\end{center}
\begin{center}\small
Table 1:\ Results of GR-PPA with fixed parameters $(\sigma_2,\sigma_3,s)=(0.2,0.2,10)$.
\end{center}
\vskip2mm
\begin{center}
\begin{tabular}{ccccccccccc}
\hline
 $s$&&  IT&&  CPU &&IER   &&OER &&CER  \\
\hline
 10 && 141  && 4.20  && 6.05e-9 &&1.60e-12  && 7.80e-9\\
 12 && 194  &&   5.49   && 2.39e-9  &&2.41e-12 & & 8.73e-9\\
 15 && 257  &&  7.34  &&  2.42e-9&&  1.53e-12& & 9.92e-9  \\
 17 &&  297  &&  8.58    &&   5.10e-9&&   2.45e-12& &  6.70e-9  \\
 20 && 363 &&  10.42 &&  6.55e-9& &2.24e-12&&  4.47e-9 \\
 22 && 435 && 12.49 &&  2.09e-9& &1.97e-12&&  8.66e-9 \\
  25 && 487 &&  14.03 &&  4.62e-9&&  7.26e-13& & 5.19e-9 \\
   27 &&510   && 15.24   &&   9.37e-9&&   1.74e-12& &  9.41e-9 \\
 30 &&  585 &&  16.93 & &  6.94e-9&& 3.57e-12&& 7.10e-9\\
  35 && 692 &&   19.80 &&  5.87e-9&&  5.75e-13&&  5.13e-9 \\
 40 && 782  &&    22.22 &&  6.54e-9&& 2.41e-12&&  9.37e-9 \\
 45 &&873  &&  25.24 &&  6.05e-9&&  2.61e-13&&  8.32e-9 \\
\hline
\end{tabular}
\end{center}
\begin{center}\small
Table 2:\ Results of GR-PPA with fixed parameters $(\sigma_1,\sigma_2,\sigma_3)=(0.178,0.178,0.178)$.
\end{center}

Secondly, we would like to use  our GR-PPA (Algorithm 2.1) with the tuned parameters   to  compare  with  two state-of-the-art
algorithms:   the Proximal Jacobian decomposition of the ALM \cite{HeXuYuan2016} (``PJALM") and
the splitting algorithm   of the ALM   \cite{HeTaoYuan2015}
(``HTY"). The default parameter $\mu=2.01$ and $H=\beta \mathbf{I}$  are used for   HTY   \cite{HeTaoYuan2015}, and
  the proximal parameter of PJALM \cite{HeXuYuan2016} is fixed as 2 as suggested
by the theory therein.
Experimental results of these   algorithms under different tolerances are reported in Table 3.
Furthermore,  with fixed   tolerances $(\epsilon_1,\epsilon_2,\epsilon_3)=(10^{-12},10^{-14},10^{-12})$,
convergence curves of the residuals  IER and OER  against the number of iterations  by different algorithms   with  different starting iterates are respectively depicted  in Figs. 1-2.
\begin{figure}[htbp]
 \begin{minipage}{1\textwidth}
 \def\figurename{\footnotesize Fig.}
 \centering
  \resizebox{12cm}{10cm}{\includegraphics{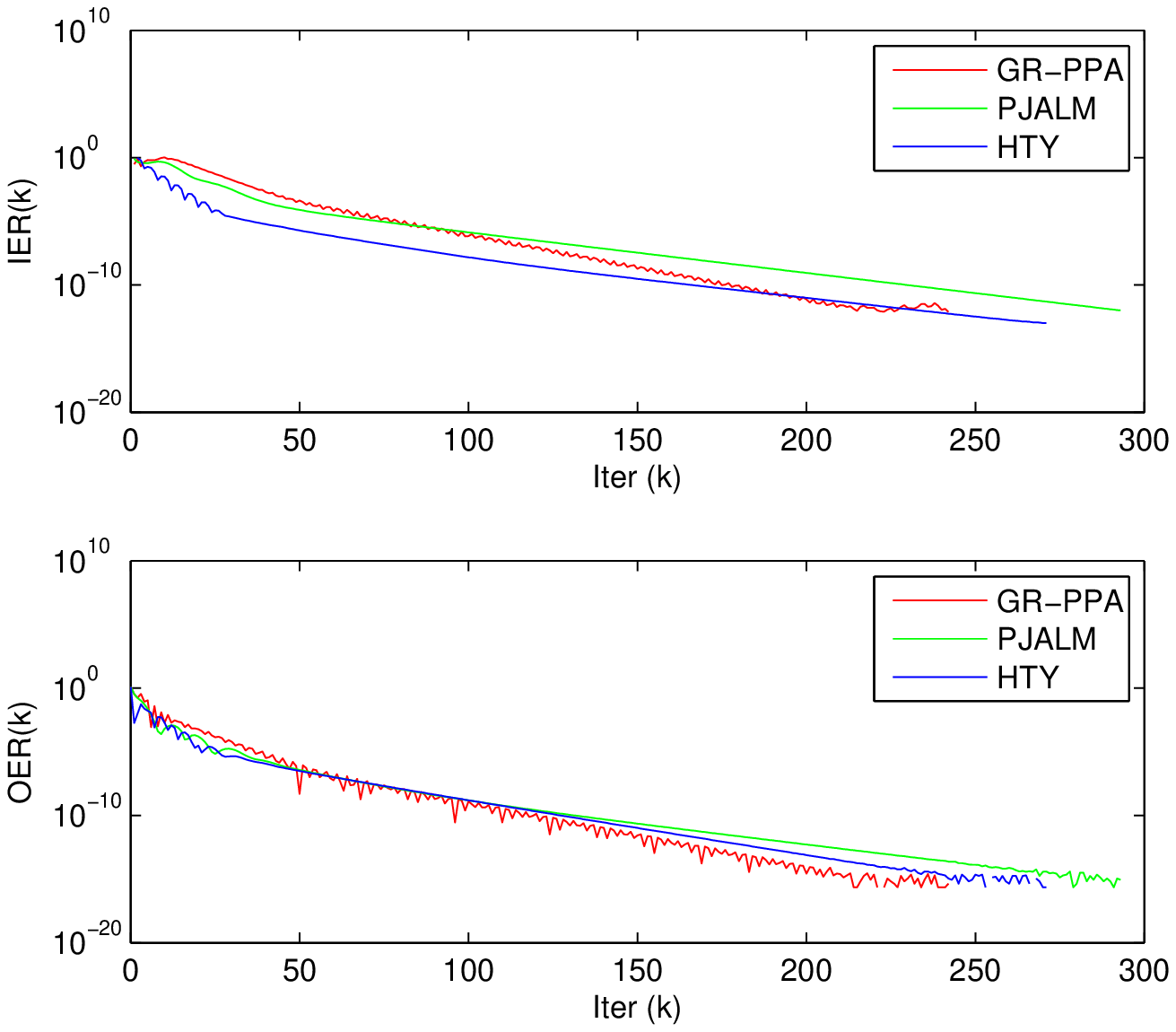}}
\caption{\footnotesize  Convergence curves of the residuals  IER and OER by different algorithms with   initial values $(X^0,S^0, L^0,\lambda^0)=(\textbf{I},4\textbf{I},3\textbf{I},\textbf{0})$.}
\resizebox{12cm}{10cm}{\includegraphics{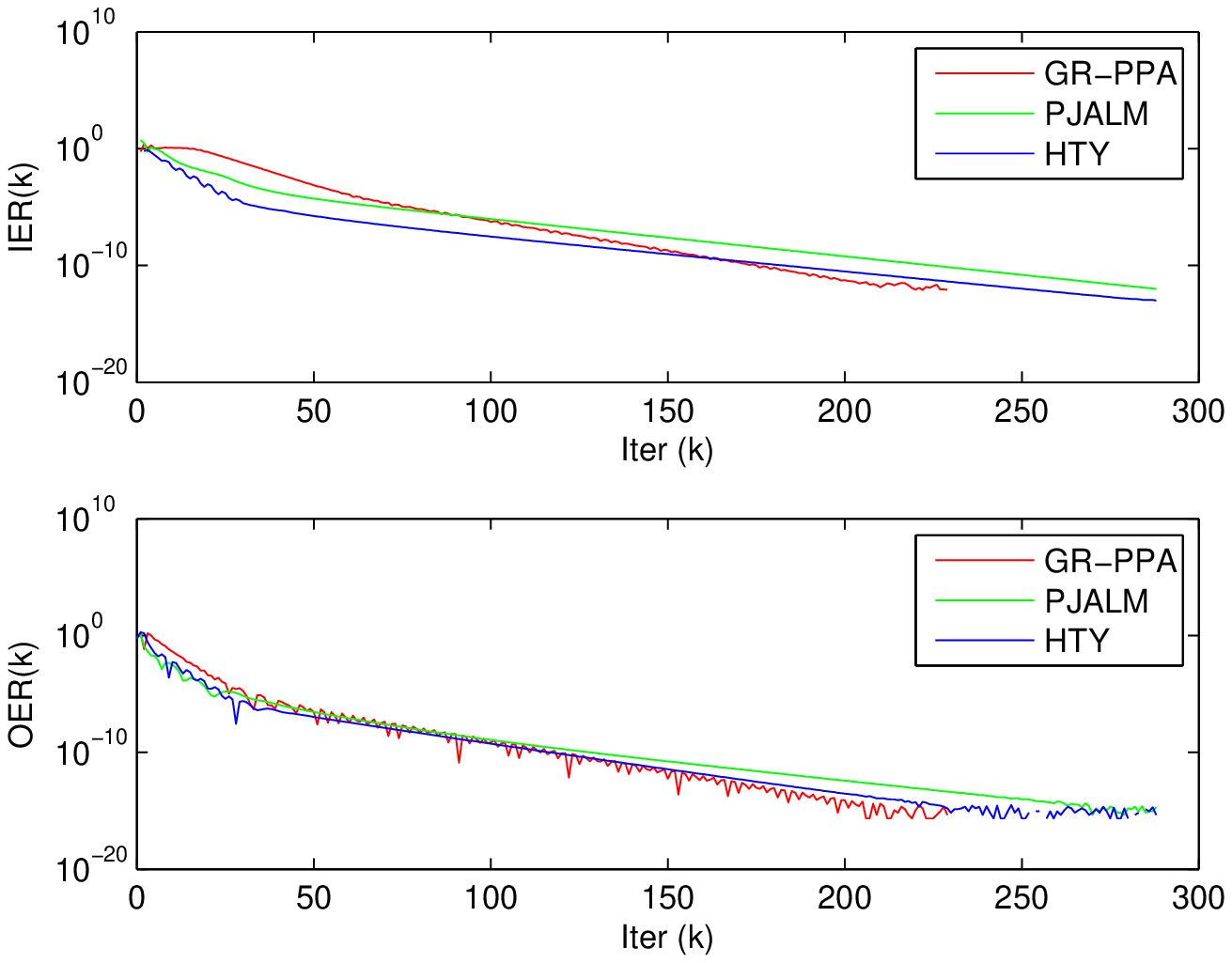}}
\caption{\footnotesize Convergence curves of the residuals  IER and OER by different algorithms with   initial values $(X^0,S^0, L^0,\lambda^0)=(2\textbf{I},4\textbf{I},2\textbf{I},\textbf{I}/2)$.}
   \end{minipage}
\end{figure}

\begin{sidewaystable}[h]
\centering % centering table
\vskip2mm\scriptsize{\begin{tabular}{ccccccccccccccccccc}
\hline\\
Parameter&\multicolumn{5}{c}{GR-PPA }& &\multicolumn{5}{c}{PJALM} & &\multicolumn{5}{c}{HTY}\\\\
\cline{2-6}\cline{8-12}\cline{14-18}{\begin{tabular}{c} $\epsilon_1(\epsilon_2=10^{-10}$\\ \quad$/\epsilon_3=10^{-6})$ \end{tabular}} & IT&  CPU &IER   &OER &CER & & IT&  CPU &IER   &OER &CER & & IT&  CPU &IER   &OER &CER\\\\
$10^{-4}$&  110  & 3.15 &  2.55e-7 & 2.89e-11 & 2.45e-7 &&
            132  & 3.55 &  1.28e-7 & 9.50e-11 & 5.59e-8 &&
            128  & 3.36 &  1.51e-9 & 9.26e-11 & 1.73e-8 \\
            \\
$10^{-7}$& 116  & 3.28 &  9.39e-8 &  7.73e-11 &  3.09e-7  &&
           136  & 3.67 &  9.53e-8 &  6.88e-11 &  4.17e-8  &&
           128  & 3.36 &  1.51e-9 &  9.26e-11 &  1.73e-8  \\
           \\
$10^{-12}$& 215  &  6.08   &  9.25e-13& 4.44e-16 &  5.49e-12  &&
            293  &  7.63   & 9.96e-13 & 8.87e-16 &  4.35e-13 &&
            234  &  6.21   & 9.56e-13 & 2.44e-15 &  1.18e-11 \\
            \\
           \hline
{\begin{tabular}{c}$\epsilon_2(\epsilon_1=10^{-6}$\\  \quad$/\epsilon_3=10^{-6})$ \end{tabular}}\\\\
$10^{-4}$& 105  & 2.72 &  4.17e-7 &  1.46e-10 &  6.58e-7  &&
           107  & 3.08 &  9.34e-7 &  9.54e-10 &  4.06e-7 &&
           79   & 2.03 &  1.03e-7 &  1.35e-8  &  9.20e-7 \\
           \\
$10^{-11}$& 124  & 3.54 &  5.11e-8  &  3.31e-12 & 5.86e-8  &&
            161  & 4.24 &  1.53e-8  &  9.76e-12 & 6.70e-9 &&
            151  & 4.01 &  2.83e-10 &  9.40e-12 & 3.41e-9 \\
            \\
$10^{-15}$& 214  & 6.15  &  1.36e-12&   4.44e-16 &  1.71e-12 &&
            279  & 7.49  &  2.74e-12&   2.22e-16 &  1.21e-12 &&
            243  & 6.48  &  5.20e-13&   8.87e-16 &  6.42e-12 \\
            \\
\hline
{\begin{tabular}{c}$\epsilon_3(\epsilon_1=10^{-6}$\\ \quad$/\epsilon_2=10^{-8})$ \end{tabular}}\\\\
$10^{-8}$ & 141    & 4.09   &6.05e-9   & 1.60e-12 & 7.80e-9 &&
            156    & 5.12   & 2.20e-8  & 1.43e-11 & 9.65e-9 &&
            136    & 3.89   & 8.32e-10 & 4.17e-11 & 9.76e-9 \\
            \\
$10^{-11}$&  200   & 6.44    & 7.24e-12 & 2.44e-15 & 6.41e-12  &&
             251   & 7.19    & 2.12e-11 & 1.26e-14 & 9.33e-12  &&
             237   & 6.96    & 7.83e-13 & 3.10e-15 & 9.66e-12  \\
            \\
$10^{-12}$& 225  & 7.31  &  1.50e-12 &  6.66e-16  &  9.87e-13 &&
            282  & 8.65  &  2.22e-12 &  3.11e-15  &  6.70e-13 &&
            271  & 8.10  &  9.94e-14 &  2.22e-16  &  9.59e-13 \\
            \\
\hline % inserts single-line
\end{tabular}
\begin{center}\small
Table 3:\ Comparative results of the problem (\ref{Sec3-prob})  by different algorithms
under different tolerances.
\end{center}}
\end{sidewaystable}

It is clear   from Table 3
that under  higher tolerances our GR-PPA  performs significantly better than the others in terms of   the number of
iterations and the CPU time, although   HTY could perform  slightly better than GR-PPA under some lower
tolerances. Besides, the comparative convergence curves depicted in Figs. 1-2  show that GR-PPA
converges  faster than the rest two  algorithms for different starting feasible points, which illustrates
that    GR-PPA could performs well for any initial iterative values. Reported results of Table 1 and
convergence curves in Figs. 1-2 demonstrate the effectiveness and robustness of the proposed algorithm.

\section{Conclusion remarks}
In this paper, we develop a relaxed parameterized PPA  for solving a class of
separable convex minimization problems.  The global convergence and a worst-case sub-linear convergence
rate of the algorithm are established. Numerical experiments on testing a sparse matrix minimization
problem in   statistical learning also verify that our proposed algorithm outperforms  two popular
algorithms  when properly choosing the relaxation factor and the  parameters in the proximal matrix.

Note that,  in Sec.3, we only take a three-block problem of (\ref{Sec1-prob}) for an example to
investigate the performance of our algorithm,  since its performance for the two-block {\it lasso} problem had been verified in \cite{BaiZhang2017} compared with the classical ADMM and the relaxed PPA \cite{GuHeYuan2014}. For the problem with more than three variables,
the parameters in  GR-PPA need to  adjust afresh via experiments. From the framework of   PPA
in (\ref{Sec2-PPA}), one may have other  choices to construct a novel parameterized proximal matrix to
develop the corresponding  PPA only if it is symmetric positive definite. Besides,   the proposed
algorithm is  applicable to the separable convex optimization problem where the coefficient matrices in the linear
constraint have full column rank. Hence,  whether there exists a PPA for
the non-separable case with nonlinear constraints? In addition, we have tried our best to analyze the algorithmic   convergence rate in the non-ergodic sense but without any results.
 These questions  need  further investigations in the future.

 \vskip 3mm \noindent{\large\bf Acknowledgements}\vskip 2mm
The authors wish to thank  the
Editor-in-Chief Prof. Choi-Hong Lai  and the  anonymous referees  for providing their valuable suggestions, which have significantly improved
the quality of our paper.

\end{document}